\documentclass[11pt,twoside,a4paper]{amsart}
\usepackage{amssymb}
\usepackage{a4wide}
\date{\today}
\usepackage[all,cmtip]{xy}

\usepackage{xcolor}

\def\e{{\epsilon}}

\def\w{\wedge}

\def\C{{\mathbb C}}
\def\w{{\wedge}}
\def\P{{\mathbb P}}

\def\M{{\mathcal M}}
\def\codim{{\rm codim\,}}

\def\K{{\mathcal K}}

\def\O{{\mathcal O}}

\def\L{{\mathcal L}}

\def\L{{\mathcal L}}
\def\U{{\mathcal U}}
\def\Pk{{\mathbb P}}
\def\Ok{{\mathcal O}}

\def\1{{\bf 1}}

\def\J{{\mathcal J}}

\def\k{{\kappa}}

\def\Mboll{{\mathring M}}

\DeclareMathOperator{\mult}{mult}
\def\m{{\mathfrak m}}

\let\Im\relax
\DeclareMathOperator{\Im}{\mathcal{I}\mathit{m}}

\newtheorem{thm}{Theorem}[section]

\newtheorem{prop}[thm]{Proposition}

\theoremstyle{definition}

\newtheorem{remark}[thm]{Remark}

\newtheorem{ex}[thm]{Example}

\numberwithin{equation}{section}

\begin{document}

\title[Analytic representation of the Buchsbaum-Rim multiplicity]
{Analytic representation of the \\ Buchsbaum-Rim multiplicity}

\date{\today}

\author{Mats Andersson \&  Richard L\"ark\"ang \& Rahim Nkunzimana \& Elizabeth Wulcan}

\address{Department of Mathematics\\Chalmers University of Technology and University of Gothenburg\\
S-412 96 G\"OTEBORG\\SWEDEN}

\email{matsa@chalmers.se}
\email{larkang@chalmers.se}
\email{rahimn@chalmers.se}
\email{wulcan@chalmers.se}



\thanks{}

\begin{abstract}
Let  $E$ and $F$ be Hermitian vector bundles over a
complex manifold $X$ and let $g\colon E\to F$
be a holomorphic morphism. In a previous work by the first author was presented a
Poincar\'e-Lelong type formula with a residue term $M^g$, defined as
the pushforward of certain generalized Monge-Amp\`ere-products on the
projectivization of $E$.
We prove that in case $g$ is injective except at the point  $x\in X$,
the current $M^g$ is precisely of the form $e[x]$, where $e$ is the
Buchsbaum-Rim multiplicity of the cokernel of the adjoint
map $g^*$, and $[x]$ denotes the Lelong current at $x$.
\end{abstract}

\maketitle

\section{Introduction}
Let $X$ be a complex manifold of dimension $n$ and let $\J_x$ be an
ideal in the local ring $\Ok_x = \Ok_{X,x}$ 
such that $\Ok_x/\J_x$ has
finite length, or equivalently, support at $\{x\}$.
Then the classical Hilbert-Samuel multiplicity $e(\J_x)$ is an important measure of the ``order of vanishing''
of $\J_x$ at $x$.
Algebraically, it may be represented as
a measure of the asymptotic
growth of the length of $\Ok_x/\J_x^k$
as $k \to \infty$. Geometrically, it may be described
as an
intersection number, cf.\ Section ~\ref{hsmult}.
If $f=(f_1,\dots,f_s)$ is a tuple of generators of $\J_x$ and
\begin{equation}\label{umgange}
  M^f= \1_{\{x\}} (dd^c \log |f|^2)^n,
  \end{equation}
where $  (dd^c \log |f|^2)^n$ is the Bedford-Taylor-Demailly
Monge-Amp\`ere product of $f$, then
a version of King's formula
states that
\begin{equation}\label{barer}
    M^f = e(\J_x)[x],
\end{equation}
where $[x]$ denotes the Lelong current at $x$, cf.\ \cite{King}
and e.g., \cite[Section III.8-9]{Dem}.

\smallskip

Assume that $\M_x$ is an $\Ok_x$-module of finite length; then $\M_x$
can be represented as $\Ok^r_x/\K_x$, where $\K_x$ is a submodule of
$\Ok_x^r$. In \cite{BR2}
Buchsbaum-Rim introduced an algebraically defined
multiplicity\footnote{In the literature there are two conventions. Sometimes $e(\M_x)$ is
referred to as the Buchsbaum-Rim multiplicity of $\K_x$ rather than
$\M_x$.} $e(\M_x)$ of $\M_x$, that generalizes the Hilbert-Samuel
multiplicity  in the sense that if $\M_x$ is of the form $\Ok_x/\J_x$, then $e(\M_x)=e(\J_x)$.
Later Kleiman-Thorup \cite{KT1, KT2} and Henry-Merle \cite{HM} gave a
geometric representation of this so-called Buchsbaum-Rim multiplicity
as an intersection number.
The aim of this note is to present an analytic representation of
$e(\M_x)$ that extends \eqref{barer}.

Let $E$ and $F$ be Hermitian holomorphic vector bundles over $X$ and
$g\colon E\to F$ a morphism that is injective outside the analytic
set $Z$.
In \cite{A1} was introduced, as a residue term in a Poincar\'e-Lelong
type formula, a current $M^g$ that generalizes
\eqref{umgange}.
It is defined in terms of generalized Monge-Ampère
products on the projectivization
$\P(E)$ of $E$.
It has support on $Z$, and is a sum
$M^g_\k+M^g_{\k+1}+\cdots +M^g_n$, where $\k$ is
the codimension of $Z$, and $M^g_k$ is closed
and has bidegree $(k,k)$.
If $E$ is a trivial line bundle and $F$ is a trivial bundle of rank
$s$, then $g$ may be identified with a tuple $f=(f_1,\ldots, f_s)$ of
holomorphic functions. If the zero-set of $f$ equals $\{x\}$,
then $M^g$ coincides
with the current \eqref{umgange}, cf.\ Examples \ref{kattvakt} and
\ref{generalisering}.

Let $g^* \colon F^*\to E^*$ denote the adjoint map of $g$. Then the sheaf
$\M:=\Ok(E^*)/\Im g^*$ has support on $Z$ and in particular, if
$Z=\{x\}$, $\M_x$ is an $\Ok_x$-module of finite length. Note
that any coherent $\O_X$-module locally can be represented like this.
There are well-defined multiplicities, or Lelong numbers, $\mult_xM^g_k$, of the
currents $M^g_k$
that
are non-negative
integers that
only depend on the sheaf
$\M$ and not the specific choice of $g$, see
\cite[Theorem~10.1]{A1}.
Our main result is the following generalization of \eqref{barer} that
gives an analytic representation of the
Buchsbaum-Rim multiplicity.
\begin{thm}\label{thm1}
  Assume that $g\colon E\to F$ is a morphism of Hermitian vector
  bundles over $X$ that is injective outside $x$. 
Let $\M_x=\Ok(E^*)_x/\Im g_x^*$ be the corresponding $\Ok_x$-module of finite length.
Then
\begin{equation}\label{evighet}
M^g=e(\M_x)[x],
\end{equation}
where $e(\M_x)$ denotes the Buchsbaum-Rim multiplicity of $\M_x$.
\end{thm}

Since $M^g$ is a global current Theorem ~\ref{thm1} yields a global
representations of Buchsbaum-Rim multiplicities. This will be further
explored in an upcoming paper.

The main ingredient in the proof of Theorem~\ref{thm1}  is a formula for $M^g$ as the
pushforward of a current on the blowup of $\P(E)$ along a certain ideal
sheaf $\J$ associated with $g$, Theorem
~\ref{tomat}.
From this we get geometric representations of $\mult_x M^g_k$
as certain intersection numbers.
When $Z=\{x\}$ this description
coincides with the intersection formula for the Buchsbaum-Rim
multiplicity due to Henry-Merle, \cite{HM}, and Kleiman-Thorup,
\cite{KT1, KT2}, and hence \eqref{evighet} follows.
If $E$ is a trivial line bundle, then $\M$ is of the form $\Ok_X/\J$
and $M^g_k$ equals $\1_Z$ times a generalized Monge-Amp\`ere
product in the sense of \cite{A3}. In \cite{ASWY} it was proved that $\mult_x
M^g_k$ in this case is the $k$th
so-called Segre number of $\J_x$.
From Theorem ~\ref{tomat} we recover this result through
the intersection formula for Segre numbers from
\cite{GG}.
In the general
case it would be interesting to compare $\mult_x M^g_k$ to the
generalized Buchsbaum-Rim multiplicities by Callejas-Bedregal and
Jorge P\'erez, \cite{CBP}.

\smallskip

In Section ~\ref{prel} we review the definition of the current $M^g$
after providing some necessary background. In Section ~\ref{multiplicera} we recall the geometric
definitions of the Hilbert-Samuel and Buchsbaum-Rim multiplicities as
well as the Segre numbers.
The core of the paper is Section ~\ref{sec:MgNew}, where we present the formula
for $M^g$ and derive expressions for its multiplicities. In particular,
we prove Theorem ~\ref{thm1}.
In Section ~\ref{vogelsektion} we construct a current $M^{g, \gamma}$,
using so-called Vogel cycles, that has the same Lelong numbers as
$M^g$ and that is well suited for computations.
In Section ~\ref{advent} we use this to
compute $M^g$ for an explicit choice of $g$.

\section{Preliminaries}\label{prel}
Throughout the paper, unless otherwise mentioned, $X$ will be a
complex manifold of dimension $n$.

\subsection{Generalized Monge-Amp\`ere products}\label{mf}
In this section we will recall the
construction of generalized Monge-Amp\`ere products from \cite{A3},
and we will allow
$X$ to be a reduced analytic space of pure dimension $n$, cf.\ \cite{ASWY}.
For a reference on currents on analytic spaces see, e.g.,
\cite[Section~2]{AESWY}.
Assume that $h \colon W'\to W$ is a proper holomorphic map between analytic
spaces, $\eta$ is a smooth form on $W$, and $\mu$ is a current on
$W'$. Then
\begin{equation}\label{pushpull}
  \eta \wedge h_*\mu = h_*( h^* \eta \wedge \mu).
\end{equation}
Moreover, if $\mu$ has order $0$ and $V$ is a subvariety of $W$, then
$h_*\mu$  has order $0$ and
\begin{equation}\label{ordobild}
  \1_Vh_* \mu = h_*(\1_{h^{-1} V} \mu).
\end{equation}

Let $f$ be a holomorphic section of a Hermitian vector bundle $F$ over
$X$
and let $\J$ be the ideal sheaf defined by $f$, with zero-set $Z$ of
codimension $\kappa$.
We recursively define generalized Monge-Amp\`ere products
\begin{equation}\label{laban}
[dd^c\log|f|^2]^{0}=1, \quad
[dd^c\log|f|^2]^{k+1}=dd^c(\log|f|^2 \1_{X\setminus
  Z}[dd^c\log|f|^2]^{k}), \ k=0, 1,2, \ldots.
\end{equation}
For $k\leq \kappa$ this coincides with the definition of the classical
Bedford-Taylor-Demailly Monge-Amp\`ere product $(dd^c \log |f|^2)^k$.

Note that if the definition \eqref{laban} makes sense it is
functorial, i.e., if $\pi \colon Y\to X$ is a modification, then
\begin{equation*}
  [dd^c\log|f|^2]^{k}=\pi_* [dd^c\log|\pi^* f|^2]^{k},
\end{equation*}
cf.\ \eqref{pushpull} and \eqref{ordobild}.
Also, note that if $f\equiv 0$, then $[dd^c \log |f|^2]^k=0$ for $k>0$;
therefore we may assume that $f\not\equiv 0$ (on each irreducible component of $X$).
Since any two modifications can be dominated by a third one, to prove that
$[dd^c\log|f|^2]^{k}$ are well-defined currents on $X$ it is enough to
prove that $[dd^c\log|\pi^* f|^2]^{k}$ are well-defined currents on $Y$
when
$\pi\colon Y\to X$ is a principalization of $\J$.
In other words, it is enough to prove that \eqref{laban} are
well-defined  when $f$ is principal, i.e., of the form
$f=f^0f'$, where $f^0$ is a section of the line bundle $\L\to X$ defining the exceptional divisor $D$
and $f'$ is a nonvanishing section of $\pi^* F\otimes \L^{-1}$. It
follows that
$0\to \L\stackrel{f'}{\to} F$
is a pointwise injective morphism, and hence $\L$ inherits a Hermitian
metric from $F$. Notice that
$|f^0|=|f^0 f'|=|f|.$
Now,  the classical Poincar\'e-Lelong formula gives
\begin{equation}\label{gamas1}
  dd^c \log |f|^2 = dd^c\log|f^0|^2 =[D]+\omega_f,
\end{equation}
where $[D]$ denotes the Lelong current of $D$, and $\omega_f$ is the
first Segre form, or equivalently minus the first Chern form, of $\L$.
In view of \eqref{gamas1},
$\omega_f= dd^c\log|f|^2$
in $X\setminus Z$.
  It follows that
\[
  [dd^c \log |f|^2]^k
 =[D]\w \omega_f^{k-1}+\omega_f^k, \quad k=1,2,\ldots,
\]
are well-defined closed $(k,k)$-currents of order $0$ and hence so are
\eqref{laban} for general $f$.

The currents $[dd^c \log |f|^2]^k$ can also be obtained as limits of
smooth forms, e.g.
\begin{equation}\label{dalmatin}
[dd^c \log |f|^2]^k=\lim_{\epsilon\to 0} \log (dd^c \log (|f|^2+\epsilon))^k.
  \end{equation}

We let
\begin{equation*}
  M^f_k=\1_Z [dd^c \log |f|^2]^k, \quad \quad M^f = \sum_{k=0}^n
  M^f_k.
\end{equation*}
Then for each $k$, $M^f_k$ is a closed $(k,k)$-current of order $0$ with support on
$Z$, and thus it vanishes if $k<\kappa$. In particular, $M^f_0=\1_Z$
vanishes unless $f\equiv 0$ on (some component of) $X$.
Also, $M^f$ is functorial, i.e., $M^f=\pi_* M^{\pi^* f}$ when $\pi$ is a
modification.
If $f$ is principal (and $f \not\equiv 0$), then
\begin{equation}\label{kaninen}
  M^{f}= 
 [D]\w \frac{1}{1-\omega_f}.
\end{equation}

\begin{remark}\label{barnvakt}
Note that if $a$ is a nonvanishing holomorphic function, then
$dd^c \log|af|^2 = dd^c \log |f|^2$, so $M^{af} = M^f$.
\end{remark}

\begin{ex}\label{kattvakt}
Assume that $f=(f_1,\ldots, f_s)$ is a tuple of germs of holomorphic
functions at $x$
generating an ideal $\J_x\subset \Ok_x$ such that $\Ok_x/\J_x$ has
finite length as in the introduction.
Then we
can identify $f$ with a holomorphic section of $F=X\times \C^s$. If we
equip $F$ with the trivial metric, then $M^f$ coincides with
the current \eqref{umgange}.
  \end{ex}


\subsection{Proper intersections}\label{solsken}
Given an analytic cycle $W=\sum a_j W_j$, where $W_j$ are the irreducible
components, we let $[W]=\sum a_j [W_j]$ be the corresponding Lelong
current. Moreover, we let $|W|=\bigcup W_j$ denote the support.

Assume that $Z_1, \ldots, Z_m$ are analytic cycles on $X$ of pure
codimensions $p_j$ that intersect properly, i.e., the intersection of
their supports $V$ have codimension $p_1+\cdots + p_m$. Then there is a
well-defined proper intersection cycle
$
  Z_m\cdots Z_1=\sum a_j V_j,
$
where $V_j$ are the irreducible components of $V$ and $a_j$ are
(integer) intersection multiplicities as defined in e.g., \cite[\S
12]{Ch},
that is
commutative in the factors.
Also there is a well-defined wedge product $[Z_m]\w\cdots\w
[Z_1]$ and this current coincides with $[Z_m\cdots Z_1]$, see, e.g.,
\cite{Ch}.

We will use that if $Z$ and $W$ are cycles that intersect properly, then
\begin{equation}\label{dimensionen}
  \1_{X\setminus |Z|}[W]=[W].
  \end{equation}

Assume that $D_1, \ldots, D_m$ are divisors such that $D_1,\ldots, D_k$
intersect
properly for each $k$, and moreover assume that they are locally defined by the holomorphic functions
$h_1, \ldots, h_m$.
Then for $k=1,\ldots, m$ the recursively defined (Bedford-Taylor-Demailly) Monge-Amp\`ere products
\begin{equation*}
  \log |h_k|^2 dd^c \log |h_{k-1}|^2\w \cdots \w dd^ c \log |h_1|^2
\end{equation*}
and
\begin{equation}\label{imponerad}
  dd^c \log |h_k|^2\w\cdots \w dd^ c \log |h_1|^2:=
  dd^c  \big (\log |h_k|^2 dd^c \log |h_{k-1}|^2\w \cdots \w dd^ c \log |h_1|^2\big )
\end{equation}
are well-defined
and the latter coincides with $ [D_k]\w\cdots \w [D_1]$.

If $D_1$ is irreducible and reduced, and $i \colon D_1\hookrightarrow X$ is the inclusion of $D_1$ in $X$ and we use the shorthand
notation $[f]$ for the divisor defined by the holomorphic function (or
section of line bundle) $f$, then
\begin{equation}\label{middag}
  \log |h_k|^2[h_{k-1}]\w\cdots \w [h_1]=
  i_*(\log |i^*h_k|^2[i^* h_{k-1}]\w\cdots \w [i^* h_2] )
\end{equation}
so that $\log |h_k|^2[h_{k-1}]\w\cdots \w [h_1]$, and analogously
$[h_{k}]\w\cdots \w [h_1]$, can be seen as the pushforward of an intrinsic current on
$D_1$.

If one of the divisors is the graph of a holomorphic function, the
intersection product may be calculated as a product of the remaining divisors restricted
to the graph. More precisely, assume that we
have local coordinates $z=(z',z_n)$ and that $h_j(z)=z_n- g(z')$.
Let $h'_k (z') = h_k(z', g(z'))$ for $k\neq j$.
Then
\begin{equation} \label{eq:MArestriction}
 [h_m] \wedge \cdots \wedge [h_1] = [h_m']\wedge \cdots \wedge
 [h_{j+1}']\wedge [h_j]\w[h'_{j-1}]\w\cdots\w [h'_1].
\end{equation}
This follows from \eqref{middag} applied to $D_j$ and the fact that
the Monge-Amp\`ere product is commutative in the factors.
For a general reference on Monge-Amp\`ere products see, e.g.,
\cite[Chapter~III]{Dem}.

\subsection{Generalized cycles and intersection of divisors}\label{genz}
The generalized Monge-Amp\`ere products \eqref{laban} and the currents
$M^f$, as well as (Lelong currents of) analytic cycles, are examples
of
generalized cycles. The $\mathbb Z$-module of generalized cycles
was introduced in \cite{AESWY}; it is the smallest class of
currents on analytic spaces that is closed under multiplication by
components of Chern forms and under direct images of proper
holomorphic maps. It is closed under multiplication by characteristic
functions 
of analytic varieties. All its elements are currents of order $0$.

A holomorphic section $f$
of a Hermitian vector bundle over $X$ induces an endomorphism of generalized cycles on $X$ by
\begin{equation}\label{brajobbat}
 \mu \mapsto  [dd^c \log |f|^2]^k\w \mu := \lim_{\epsilon \to 0} (dd^c \log
  (|f|^2+\epsilon))^k\w \mu
\end{equation}
see Theorem~5.2 in \cite{AESWY},
cf.\ \eqref{dalmatin}.
Hence also
\begin{equation*}
 \mu\mapsto M^f_k \wedge \mu := \1_Z [dd^c \log |f|^2]^k\w \mu
\end{equation*}
is such an endomorphism. Note that $M^f_k\w \mu$ can be nontrivial
even if $M^f_k=M_k^f\w 1=0$.
In particular, if $W$ is an analytic
cycle and $D$ is a divisor on $X$ locally defined by the holomorphic
function $f$, we define
\begin{equation}\label{glomska}
  [D]\w [W]= M^f_1\w [W],
\end{equation}
cf.\ \eqref{kaninen},
and thus recursively intersections of divisors $[D_m]\w\cdots\w
[D_1]$.
It follows from \eqref{brajobbat} that
\begin{equation}\label{renovering}
[D]\w [W] = [D]\w \1_{X\setminus |D|}
[W],
\end{equation}
where the right hand side is a proper intersection; in particular,
\eqref{glomska} is independent of the choice of $f$.
Note that, if $D$ and $W$ intersect properly, then $[D]\w
[W]$ defined by \eqref{glomska} coincides with the regular proper
intersection, cf.\
\eqref{dimensionen}.

Assume that $h_1,\ldots, h_m$ are holomorphic functions
(or sections of line bundles). Then, it follows from
\eqref{brajobbat} that,
for $k<\ell$, $c\in \C$,
\begin{equation}\label{lussekatt}
  [h_\ell + c h_k]\wedge [h_{\ell -1}]\wedge \cdots \wedge [h_1]=
  [h_\ell]\wedge [h_{\ell -1}]\wedge \cdots \wedge [h_1].
  \end{equation}

If $h \colon X'\to X$ is a proper holomorphic map and $\tau$ is a generalized
cycle on $X'$, then $\mu:=h_* \tau$ is a generalized cycle and
\begin{equation}\label{lurig}
  M^f\w\mu=h_* (M^{h^*f}\w \tau).
  \end{equation}

Assume that $\mu$ is a generalized cycle. Then there is a well-defined multiplicity at $x\in X$.
Let $\xi$ be a section of a vector bundle in a neighborhood $\U
\subset X$ of $x$ such that $\xi$ defines the maximal ideal $\m$ at
$x$.
Then
\begin{equation}\label{kalsong}
\mult_x \mu := \int_\U M^\xi \wedge \mu
\end{equation}
is an integer that only depends on $\mu$ and not the choice
of $\xi$.
If $\mu$ is a positive closed current $\mult_x \mu$ coincides with the
Lelong number at $x$.

\subsection{Vogel cycles}\label{tidigast}

Assume that $f=(f_1,\ldots, f_s)$ is a tuple of holomorphic functions
in a neighborhood $\U$ of $x\in X$
(or a section of a trivial vector bundle $F=X\times \C^s$) generating a
local ideal $\J_x$ with zero-set $Z$ of codimension
$\kappa$.
For $\beta=[\beta_1:\ldots :\beta_s]\in \P^{s-1}$, let
\[
  \beta\cdot f=\beta_1f_1+\cdots +\beta_s f_s,
\]
and for $\beta=(\beta_1, \ldots, \beta_n)\in (\P^{s-1})^n$,
let
\[
 M^{f,\beta}_0=\1_Z, \quad \quad  M^{f,\beta}_k=\1_Z[\beta_k\cdot f]\wedge\cdots \wedge
  [\beta_1\cdot f], \quad k=1,\ldots, n,
\]
where the product is defined recursively as in Section ~\ref{genz},
cf.\ \eqref{renovering}.
If $\beta\in (\P^{s-1})^n$ is generic, then
$h=(h_1,\ldots, h_n)$, where $h_\ell=\beta_\ell\cdot f$, is a
so-called \emph{Vogel sequence} of $\J$ at $x$, i.e.,
after possibly shrinking $\U$,
\[
  \codim \big ( (\U\setminus Z)\cap \{h_1=0\}\cap \cdots \cap
  \{h_k=0\}\big ) = k \text{ or }\infty, \text{ for } k=1,\ldots, n.
  \]
 It follows in view of \eqref{renovering} and \eqref{dimensionen}
  that
    $M^{f,\beta}:=\sum_{k=0}^n M^{f,\beta}_k$
        coincides with the corresponding
\emph{Vogel cycle} $V^h=V^h_\kappa+\cdots +V^h_n$, where
\[
   V^h_k=
  \1_Z[h_k]\1_{X\setminus Z}
 [h_{k-1}]\w\cdots\w\1_{X\setminus Z} [h_{1}].
\]
Theorem~6.1 in \cite{ASWY} asserts that for a generic $\beta\in (\P^{s-1})^n$
\begin{equation}\label{burra}
  \mult_x M_k^{f,\beta} = \mult_x M^f_k.
  \end{equation}
  In particular,
  if $Z=\{x\}$, then
  \begin{equation}\label{hurra}
    M^{f,\beta}= (\mult_x M^f_n) [x]  = M^f.
  \end{equation}
In \cite{ASWY} this result was used to prove that $\mult_xM^f_k$ is
the $k$th Segre number
of $\J_x$, and also for explicit computations.

\subsection{The projectivization of $E$}\label{hundvakt}
Assume that $E$ is a vector bundle of rank $r$ over $X$. Let $\Pk(E)$ be the projectivization of $E$, where
the fiber at $x$ is $\Pk(E_x)$, i.e., the set of all lines in $E_x$
through the origin\footnote{There is another common convention in the
  literature where
the fiber instead consists of all hyperplanes through the origin.}, and let $p \colon
\Pk(E) \to X$ be the natural submersion.
In a local trivialization where $E=X\times \C^r$, note that $\P(E)=X\times
\P^{r-1}$. Throughout the paper we let
\[
  N=\dim \P(E)=n+r-1.
  \]
Consider the pullback $p^*E$ of $E$ to $\Pk(E)$ and let $L_E$ denote
the tautological
sub(line)bundle whose fiber at $(x,[\alpha])$ is precisely the line through
the origin of $E_x$ determined by $[\alpha]$.

If $E$ is an Hermitian vector bundle, i.e., equipped with a smooth Hermitian metric,
then $p^*E$ inherits a smooth metric, and thus so does $L_E$.
Note that for $\alpha \in E_x$,
$dd^c \log |\alpha|^2$ only depends on $[\alpha]$. Thus $dd^c\log|\alpha|^2$ is well-defined form om $\Pk(E)$.
One can verify that, in fact,
\begin{equation*}
\omega_\alpha:=dd^c\log|\alpha|^2,
\end{equation*}
is the first Segre form, or equivalently minus the first Chern form, of the induced metric on $L_E$, see, e.g., \cite[Section~4]{A1}.

\subsection{The current $M^g$} \label{sec:MgDef}Let us recall the
construction of the currents $M^g$ in \cite[Section~4]{A1}.
Assume that $E$ and $F$ are Hermitian vector bundles over $X$ of rank $r$ and $s$,
respectively, and that $g \colon E\to F$ is a holomorphic morphism; let $Z$
denote the set where $g$ is not injective.
Let $\P(E)$ be the projectivization of $E$ with tautological line
bundle $L_E$.
Note that a local section $\sigma$ of $L_E$ is of the form
\begin{equation*}
  \sigma = s(x, \alpha) \alpha
\end{equation*}
at $(x,[\alpha])\in \P(E)$, where $s(x,\alpha)$ is a holomorphic
function on $E\setminus \{\mathbf{0}\}$
that is $-1$-homogeneous in $\alpha\in
E_x\setminus\{\mathbf{0}\}$; here $\{\mathbf{0}\}$ denotes the zero section.
Therefore, we may identify sections of $L_E$ with such functions
$s(x,\alpha)$ and $\alpha$ can be considered as a section of
$p^*E\otimes L_E^{-1}$. Consequently
\[
  G(x, [\alpha]):= g\alpha
\]
can be seen as a section of $p^*F\otimes L_E^{-1}$; here $g$ denotes also the
induced morphism $p^*E\to p^*F$.
Let $a$ be a local nonvanishing section of $L_E$ defined on $U_a
\subseteq \Pk(E)$. Then on $U_a$, let
\begin{equation*}
\Mboll^{G}=M^{Ga}.
\end{equation*}
In view of Remark ~\ref{barnvakt},
$M^{Ga}$ is independent of the choice of $a$ and
thus $\Mboll^G$ is a well-defined global current on $\Pk(E)$.
Note that $\Mboll^G$ depends on the Hermitian metric on $F$, but that the metric on $L_E$ is suppressed.

Next, we let $\omega_\alpha$ be the first Segre
form of $L_E$ as in Section ~\ref{hundvakt}.
Then $M^g$ is defined as
\begin{equation} \label{Mgdef}
M^g=p_*\Big(\Mboll^G\w \frac{1}{1-\omega_\alpha}\Big),
\end{equation}
cf.\ \cite[Definition~4.1]{A1}.
It follows that  $M^g =\sum_{k=0}^n M^g_k$, where $M_k^g$ is a
closed current of order $0$ and bidegree $(k,k)$ with support on
$Z$. In particular, if $Z=\{x\}$, then
\begin{equation} \label{MgPoint}
    M^g = M^g_n = \lambda [x]
\end{equation}
for some $\lambda \in \mathbb{C}$.
Moreover, the $M^g_k$ are generalized cycles, so there are well-defined
multiplicities $\mult_xM^g_k$; in fact the multiplicities, and
in particular $\lambda$, are
non-negative integers, see
\cite[Theorem~4.4]{A1}.

\begin{ex}\label{generalisering}
  Assume that we are in the setting of Section
  ~\ref{mf}, i.e., that $f$ is a holomorphic section of a Hermitian vector
  bundle $F$ over $X$.
 Then we can regard $f$ as a morphism $g \colon E\to F$, where $E$ is a
 trivial line bundle.
 Now $\P(E)=X$ and $L_E=E$ and $Z$ is the zero-set of the ideal sheaf generated by $f$. Notice that $\Mboll^G$ coincides with the
 current $M^f$ in Section ~\ref{mf}.
 Moreover if we equip $E$ with the trivial metric, then
 $\omega_\alpha=0$ and
 \begin{equation}\label{eftermiddag}
   M^g=\Mboll^G=M^f.
   \end{equation}
  \end{ex}

  \begin{ex}\label{triviala}
Assume that $g\equiv 0$. Then $G \equiv 0$ and $\Mboll^G=\1_{\P(E)}$
so that $M^g=p_*(1/(1-\omega_\alpha))$ is smooth. In fact, $p_*(1/(1-\omega_\alpha))$ by definition
equals the Segre form $s(E)$ of $E$.
    \end{ex}

\section{Multiplicities and intersection formulas}\label{multiplicera}
In \cite{BR2} Buchsbaum and Rim generalized the classical
Hilbert-Samuel multiplicity of ideals to modules.
Given an $\Ok_x$-module $\M_x$ of finite length, or equivalently, with
support at $\{x\}$, they
defined a multiplicity $e(\M_x)$, now called the Buchsbaum-Rim multiplicity,
as the leading
coefficient of a certain generalized Hilbert-Samuel polynomial; in
fact, they considered modules over any Noetherian local
ring.
If $\M_x$ is of the form $\O_x/\J_x$, then $e(\M_x)$ coincides with
the Hilbert-Samuel multiplicity of $\J_x$.
In the more recent literature the Buchsbaum-Rim multiplicity is often
associated to a submodule $\K_x\subset \Ok_x^s$ rather than
 $\M_x=\Ok_x^s/\K_x$.
In the 90's Kleiman-Thorup \cite{KT1, KT2} and Henry-Merle \cite{HM}
independently developed a geometric theory for Buchsbaum-Rim
multiplicities; in particular, they proved an intersection formula.

Also, in the 90's Tworzewski, \cite{Twor}, Gaffney-Gassler,
\cite{GG}, and Achilles-Manaresi, \cite{AM} independently introduced
generalized multiplicities of ideals $\J_x\subset \Ok_x$ such that
$\Ok_x/\J_x$ does not necessarily have finite length. Assume that $\J_x$ has zero-set $Z$ of codimension
$\kappa$. Then there is an associated list of numbers $e_k(\J_x)$,
$k=\k,\dots,n$.  If $Z=\{x\}$ then the only nonvanishing number
$e_n(\J_x)$ coincides with the Hilbert-Samuel multiplicity $e(\J_x)$.
The definition of these so-called Segre numbers\footnote{We
  follow \cite{GG} where $e_k(\J_x)$ are called Segre numbers. There
  are also other names used; e.g., in \cite{CBP} they are called the
  Achilles-Manaresi multiplicity sequence.} in \cite{AM} was purely
algebraic, whereas in \cite{Twor}
and \cite{GG} was given a geometric definition in terms of
so-called Vogel cycles, cf.\ Section ~\ref{tidigast}. In
\cite[Section~3]{GG} there were also intersection formulas for the
Segre numbers.

We will recall the
geometric intersection formulas for these multiplicities. They are all
formulated in terms of blowups along ideals.

\subsection{The blowup along an ideal (sheaf)}
Assume that $f$ is a holomorphic section of a vector
bundle $F$ over $X$
that is not identically equal to
$0$. Then the generated ideal sheaf $\J$
is nontrivial and the
zero-set $Z$
is a proper subvariety of $X$.
Let
$\pi \colon  Y\to X$ be the (reduced) blowup of $X$ along $\J$. Then $Y$ is the
closure in $\Pk(F)$ of the graph
\begin{equation}\label{blowup}
A=\{ (x, [f(x)] ); x\in X\setminus Z\}
\end{equation}
over $X\setminus Z$, and $\pi$ is the restriction of the natural
projection $\Pk(F)\to X$, i.e., we have a factorization
\begin{equation}\label{trekant}
\xymatrix{
  Y \ar@{^{(}->}[r] \ar[rd]^{\pi} & \P(F) \ar[d]  \\
  & X,
}
\end{equation}
where the horizontal map is the natural inclusion and the vertical map
is the natural projection.

\subsection{The Hilbert-Samuel multiplicity}\label{hsmult}
Let us start by recalling the classical intersection formula for the
Hilbert-Samuel multiplicity.
Assume that $\J_x$ is an ideal in the local ring $\Ok_x=\O_{X,x}$; we
may assume that $\J_x$ is generated by a germ of a
holomorphic section of $F=X\times \C^s$.
Let $\pi \colon  Y\to X$ be the blowup along $\J$ in (some neighborhood of $x$
in) $X$ and let $D$ denote the corresponding exceptional
divisor.
Moreover, let $h$ denote the first Segre class of the tautological line bundle of the
blowup $Y\subset \P(F)=X\times \P^{s-1}$.
If $\J_x$ has
zero-set $\{x\}$, then the Hilbert-Samuel multiplicity of $\J_x$ equals
\begin{equation}\label{ursprung}
 e(\J_x) = \int_{Y} h^{n-1} D.
  \end{equation}

  \subsection{Segre numbers}\label{segremult}
Now assume that $\J_x\subset \Ok_x$ is an ideal whose zero-set $Z$
may have positive dimension.
Let us consider the
setting in the previous section.
Moreover, let $\m$ be
the maximal ideal in $\Ok_x$; we may assume that it is generated by a germ of a holomorphic section of $X\times \C^n$.
Let
$\pi_\m \colon Y_\m \to X$ be the blowup along $\m$ with exceptional divisor
$D_\m$, and let $h_\m$ denote the first Segre class of the tautological line bundle of the
blowup $Y\subset X\times \P^{n-1}$.
Next, let $\widetilde Y$ be the
blowup of $X$ along $\m \J$ with natural projections
\begin{equation}\label{diagram}
\xymatrix{
\widetilde Y \ar[d]^{\varpi_\m} \ar[r]^{\varpi} & Y_\m  \ar[d]^{\pi_\m}\\
Y \ar[r]^{\pi} & X,
}
\end{equation}
and let $\widetilde D$, $\widetilde D_\m$, $\widetilde h$, and
$\widetilde h_\m$ be the pullbacks to $Y$.
Note that we can decompose $D = D^{\{x\}} + D^{X\setminus \{x\}}$,
where each irreducible component of $D^{\{x\}}$
is mapped onto $\{x\}$ and no irreducible component of  $D^{X\setminus \{x\}}$ is
mapped onto $\{x\}$. Note that
\begin{equation*}
  [D^{\{x\}}]=\1_{\pi^{-1}\{x\}}[D], \quad\quad
  [D^{X\setminus\{x\}}]=\1_{Y\setminus \pi^{-1}\{x\}}[D],
\end{equation*}
cf.\ \eqref{dimensionen}.
Similarly we can decompose $\widetilde D$.
The intersection formulas in \cite[(3.2)]{GG} assert that
\begin{equation}\label{stackmoln}
  e_n(\J_x) = \int_{Y} h^{n-1} D^{\{x\}}, \quad \quad \quad
   e_k(\J_x) = \int_{\widetilde Y} \widetilde h_\m^{n-k-1}\widetilde
   h^{k-1}  \widetilde D_\m\cdot \widetilde D^{X\setminus \{x\}}, \quad k<n.
 \end{equation}
If $Z=\{x\}$, then $D^{\{x\}}=D$, so the formula for $e_n(\J_x)$
coincides with
\eqref{ursprung}.

\begin{remark}\label{mintgron}
We can interpret the integrals in \eqref{stackmoln} and
\eqref{ursprung} analytically in the following way:
We replace each power of a Segre class by the same power of a
representing smooth Segre form, and we think of the cycles
as Lelong currents. The integrals we then get are well-defined and provide the expected numbers.
\end{remark}

\subsection{Buchsbaum-Rim multiplicities}\label{brmult}
Assume that $\M_x$ is an $\Ok_x$-module of finite length.
We may
assume that $\M_x=\Ok(E)_x/\Im
g^*_x$, where $g \colon E\to F$ is a germ of a morphism of vector bundles
$E=X\times \C^r$ and $F=X\times
\C^s$ as in Section ~\ref{sec:MgDef}; let us use the notation from there.
Let $\J$ be the ideal sheaf defined by $G$ and let $\pi_E \colon Y\to \P(E)$
 be the blowup of $\P(E)$ along $\J$ with exceptional divisor
 $D$. Then $Y$ has dimension $N$ and thus $D$ is a subspace of dimension $N-1$ of
  \begin{equation*}
   \P(p^*F) \cong
   \P(E)\times_X \P(F) \cong X\times \P^{r-1}\times
   \P^{s-1},
 \end{equation*}
 where $\times_X$ denotes the fiber product over $X$.
For $\ell=1,2,\ldots$ let $H_E^\ell$ and $H_F^\ell$ denote the intersection of $\ell$ generic hyperplanes in $\P^{r-1}$ and
$\P^{s-1}$, respectively.
Since $\M_x$ has support at $\{x\}$, it follows that $D \subset \{ x \} \times \P^{r-1} \times \P^{s-1}$, so we may consider
$D$ as a projective algebraic cycle.
Now \cite[Proposition~4.2.3]{HM} asserts that\footnote{In \cite{HM} the intersection product in \eqref{HM-BR} is due to Snapper, \cite{Snapper}.
On a projective manifold however, this intersection product coincides with
the standard one in, e.g., \cite{Fulton}, which
follows, e.g., from Kleiman's axiomatic approach to intersection
products, \cite{Kleiman}.}
\begin{equation} \label{HM-BR}
e(\mathcal{M}_x) = \sum_{j=0}^{N-1} D \cdot H_E^j \cdot H_F^{N-1-j},
\end{equation}
cf.\ \cite{KT1, KT2}. Here the right hand side should be interpreted as
the multiplicity of the proper intersection cycles.

Let $h_E$ and $h_F$ be the first Segre classes of the tautological
line bundles of $\P(E)=X\times \P^{r-1}$ and $\P(F)=X\times \P^{s-1}$, respectively. Then, restricted to
the fiber $\{x\}\times \P^{r-1}\times \P^{s-1}$ over $x$ in
$\P(E)\times_X\P(F)$, $h_E$ is cohomologous to $[H_E]$ and $h_F$ is
cohomologous to $[H_F]$. Thus we can write \eqref{HM-BR} as
\begin{equation}\label{dimma}
  e(\M_x) = \int_{Y} h_E^{j}h_F^{N-1-j}D,
 \end{equation}
 cf.\ Remark \ref{mintgron}.

Assume that we are in the setting of Example ~\ref{generalisering}, so
that $E$ is a trivial line bundle and $\M_x=\Ok_x/\J_x$.
Then $\P(E)=X$, $Y$ is the blowup along $\J$, $h_E=0$, and thus
the right hand side of \eqref{dimma} coincides with the right hand
side of \eqref{ursprung}. Hence we get back the
%
%
intersection formula \eqref{ursprung} for $e(\M_x)=e(\J_x)$.

\section{A new formula for $M^g$} \label{sec:MgNew}
Assume that we are in the setting of Section ~\ref{sec:MgDef}, and
assume $g$ is nontrivial, i.e., $g\not\equiv 0$.
Let $\J$ be the non-trivial ideal sheaf defined by $G$ with zero-set
$Z\subset \P(E)$, and let
$\pi_E \colon Y\to \P(E)$ be the blowup of $\P(E)$ along $\J$ with
exceptional divisor $D$.
Assume that $a$ is a local nonvanishing section of $L_E$ in
$U_a\subset \P(E)$. Then $Ga$ is a local section of $p^*F\to \P(E)$
defining $\J$ and by Section ~\ref{mf}, on $U_a$
\begin{equation}\label{mittwoch}
  \Mboll^G= M^{Ga} =(\pi_E)_* \Big ([D]\wedge \frac{1}{1-\omega_{Ga}}\Big ),
  \end{equation}
  where $\omega_{Ga}$ is a smooth form on $\pi_E^{-1} U_a$ and
  \begin{equation}\label{freitag}
  \omega_{Ga}=\pi_E^* dd^c \log |Ga|^2
  \end{equation}
  in $\pi_E^{-1}(U_a\setminus Z)=\pi_E^{-1} U_a\setminus |D|$.

Over $U_a$, $Y$ is the closure of
\begin{equation}\label{orm}
A=\big \{ \big (x,[\alpha], [Ga(x, [\alpha])]\big );  \ (x,[\alpha])\in U_a\setminus
Z \big \}
\end{equation}
in $\P(p^*F)\cong \P(E)\times_X \P(F)$ and $\pi_E$ is the restriction of the natural
projection
$\P(E)\times_X \P(F) \to  \P(E)$.

Let $\pi_F \colon  Y\to \P(F)$ be the restriction of the natural
projection
$\P(E)\times_X \P(F) \to  \P(F)$, so that we have the commutative
diagram
\begin{equation}\label{diagram2}
\xymatrix{
Y \ar[d]^{\pi_F} \ar[r]^{\pi_E} & \P(E)  \ar[d]^{p}\\
\P(F) \ar[r] & X,
}
\end{equation}
let $\beta$ denote vectors in $F_x$,
and let $\omega_\beta=dd^c \log
|\beta|^2$ be the first Segre form of the tautological line bundle $L_F$ of
$\P(F)$, as
defined in Section ~\ref{hundvakt}. In view of \eqref{orm},
$$\omega_{Ga}=\pi_F^* \omega_\beta$$
in
$\pi_E^{-1}U_a\setminus |D|$.
Since $\omega_{Ga}$ and $\pi_F^*\omega_\beta$ are smooth forms on
$\pi_E^{-1}U_a$
that coincide outside a subvariety of positive
codimension, we conclude that they indeed coincide on $\pi_E^{-1}U_a$.
Thus on $\pi_E^{-1}U_a$
\begin{equation}\label{dokumentar}
\Mboll^G=(\pi_E)_*\left([D]\w\frac{1}{1-\pi_F^*\omega_\beta}\right).
\end{equation}
Since the left and right hand sides are independent of $a$,
\eqref{dokumentar} holds globally on $\P(E)$.
Hence, in view of \eqref{Mgdef}, cf.\ \eqref{pushpull},
we get the following formula for
$M^g$.

\begin{thm}\label{tomat}
Assume that $E,F$ are Hermitian vector bundles over $X$ and that $g\colon E\to
F$ is a nontrivial holomorphic morphism.
Let $p \colon \P(E)_{(x,[\alpha])}\to X$ be the projectivization of $E$, let $\J$ be the ideal sheaf on $\P(E)$ defined by
$p^*g \alpha$, let $Y\subset \Pk(E)\times_X \Pk(F)$ be the blowup
of $\P(E)$ along $\J$ with exceptional divisor $D$, let
$\pi_E \colon Y \to \Pk(E)$ and $\pi_F\colon Y\to \Pk(F)$ be the
restrictions of the natural projections from $\P(E)\times_X\P(F)$, and
moreover, let $\omega_\alpha$
and $\omega_\beta$ be the first Segre forms of the tautological line
bundles of $\P(E)$ and $\P(F)$, respectively.
Then
\begin{equation*}
  M^g=p_*(\pi_E)_*\Big([D]
  \w\frac{1}{1-\pi_F^*\omega_\beta}
  \w\frac{1}{1-\pi_E^*\omega_\alpha}\Big).
\end{equation*}
\end{thm}

Combining Theorem ~\ref{tomat} and the geometric description of the
Buchsbaum-Rim multiplicity in Section ~\ref{brmult} we get a proof of
Theorem ~\ref{thm1}.

\begin{proof}[Proof of Theorem ~\ref{thm1}]
By \eqref{MgPoint}, after possibly shrinking $X$, $M^g=\lambda
[x]$. Also, we may assume that $E$ and $F$ are trivial so that we are
in the setting of Section \ref{brmult}.
Thus, in view of Theorem ~\ref{tomat} and \eqref{dimma}
\begin{equation*}
  \lambda =\int_X M^g =
  \int_Y
  [D] \w\frac{1}{1-\pi_F^*\omega_\beta}\w\frac{1}{1-\pi_E^*\omega_\alpha}
  =
  \int_Y\sum_{j=0}^{N-1} [D]
  \wedge
  (\pi_E^*\omega_\alpha)^{j}\wedge (\pi_F^*\omega_\beta)^{N-1-j} =
  e(\M_x),
\end{equation*}
cf.\ Remark \ref{mintgron}.
\end{proof}

\begin{remark}\label{borste}
  Note that \eqref{mittwoch} holds for any principalization
  $\pi_E \colon Y\to \P(E)$ of $\J$.
  Since two nonvanishing sections of $L_E$ differ by multiplication
  by a nonvanishing holomorphic function, it follows in view
  \eqref{freitag} that $\omega_{Ga}$ defines a global smooth form
  $\omega_G$ on
  $Y$.
Now
\begin{equation*}
  M^g=p_*(\pi_E)_*\Big([D]
  \w\frac{1}{1-\omega_G}
  \w\frac{1}{1-\pi_E^*\omega_\alpha}\Big).
\end{equation*}
\end{remark}

\begin{ex}\label{blapenna}
Let us revisit Example ~\ref{generalisering},
i.e., assume that $f$
is a holomorphic section of a Hermitian vector
  bundle $F\to X$ regarded as
    a morphism $g \colon E\to F$, where $E$ is a
    trivial line bundle equipped with the trivial metric.
    Moreover assume that $f\not\equiv 0$.
Then $\J$ is the nontrivial ideal sheaf in $\P(E)=X$ defined by $f$ and
$\pi=\pi_E \colon Y\to X$ is the blowup of $X$ along $\J$,
 here $Y$ is the closure of \eqref{blowup} in $\P(F)$, and
$\pi_F$ is just the inclusion $i\colon Y \hookrightarrow \P(F)$.
Now in view of \eqref{eftermiddag} and \eqref{dokumentar},
\begin{equation*}
  M^f=\Mboll^G= \pi_*\left([D]\w\frac{1}{1-i^*\omega_\beta}\right).
\end{equation*}
   \end{ex}

   \subsection{Multiplicities of the currents $M^g_k$}\label{hundra}
As mentioned in the introduction, by \cite[Theorem ~10.1]{A1}, the
multiplicities
$\mult_x M^g_k$ only depend on the module $\M_x=\Ok(E)_x/\Im g_x$ and
not on the specific choice of $g$. Let us denote these numbers by
$e_k(\M_x)$. In
this section we derive intersection formulas for $e_k(\M_x)$
similar to the ones in Section ~\ref{multiplicera} in the
general case when $\M_x$ does not necessarily has finite length.

Recall from Example \ref{triviala} that if $g\equiv
0$, then $M^g$ is smooth so that all multiplicities are $0$. We may
therefore assume that $g$ is nontrivial.
Let $\xi$ be a section of a vector bundle $E_\m$ in a neighborhood $\U
\subset X$ of $x$ such that $\xi$ defines the maximal ideal $\m$ at
$x$.
Then, by \eqref{kalsong},
\begin{equation}\label{ballong}
\mult_x M^g_k = \int_\U M^\xi \wedge M^g_k=\int_\U M^\xi_{n-k} \wedge
M^g_k.
  \end{equation}

As in Sections ~\ref{brmult} and ~\ref{sec:MgNew}, let $\J$ be the
ideal sheaf defined by $G$ and let $\pi_E\colon Y  \to \P(E)$ be the
blowup of $\P(E)$ along $\J$ with exceptional divisor $D$.
Moreover, let $\pi\colon Y_{p^*\m}\to \P(E)$ be the blowup of
$\P(E)$ along $p^*\m$ with exceptional divisor $D_\m$, and let
$\widetilde Y$ be the blowup of $\P(E)$ along $\J p^*\m$.
Then we have the commutative diagram
\begin{equation}\label{stortdiagram}
\xymatrix{
& \widetilde Y \ar[ld]_{\varpi} \ar[rd]^{\varpi_\m} & \\
 Y \ar[d]_{\pi_F} \ar[rd]^{\pi_E} & & Y_{p^*\m}\ar[ld]_{\pi} \ar[d]^{\pi_{E_\m}}  \\
 \P(F) \ar[rd] & \P(E) \ar[d]_{p} & \P(E_\m)
 \ar[ld]^{p_\m}
 \\
 & X,
  }
\end{equation}
where the maps are (restrictions of) the natural
projections, cf.\ \eqref{diagram} and \eqref{diagram2}.

Let $\omega_\gamma$ be the first Segre form of the tautological line
bundle $L_{E_\m}$ of $\P(E_\m)$, as defined in Section
~\ref{hundvakt}.
Then
\begin{equation}\label{inflammation}
  M^{\pi^*p^*\xi}=
   [D_\m]\w\frac{1}{1-\pi_{E_\m}^*\omega_\gamma}.
 \end{equation}
To see this, let $\pi_\m\colon Y_\m \to X$ be the blowup of $X$ along $\m$
with exceptional divisor $\widehat D_\m$
and consider the following commutative diagram that extends the lower right square in
\eqref{stortdiagram}
\begin{equation*}
\xymatrix{
  Y_{p^*\m} \ar@{^{(}->}[r] \ar[rd]^{\pi_{E_\m}} \ar[d]_{\hat p} &
  \P(E)\times_X \P(E_\m) \ar[r]
  \ar[d] &
  \P(E) \ar[d]^{p} \\
Y_\m \ar@{^{(}->}[r]^{i} & \P(E_\m) \ar[r]^{p_\m} & X.
  }
\end{equation*}
Here the horisontal and vertical maps are the natural inclusions and
projections, cf.\ \eqref{trekant}.
Then in view of Section \ref{mf} and the beginning of Section
\ref{sec:MgNew}, cf.\ Example
~\ref{blapenna},
\begin{equation*}
   M^{\pi_\m^*\xi}=
 [\widehat D_\m]\w\frac{1}{1-i^*\omega_\gamma},
 \end{equation*}
  and now \eqref{inflammation} follows,
 \begin{equation*}
   M^{\pi^*p^*\xi}=M^{\hat p^* \pi_\m^*\xi}=
 [D_\m]\w\frac{1}{1-\hat p^* i^*\omega_\gamma}
 =
 [D_\m]\w\frac{1}{1-\pi_{E_\m}^*\omega_\gamma}.
 \end{equation*}

Let $\widetilde \U$ denote the preimage of $\U$ in $\widetilde Y$ and let $\widetilde D$ and $\widetilde D_\m$ denote the pullbacks of $D$ and
$D_\m$, respectively, to $\widetilde Y$.
Then it follows from  \eqref{lurig} and Theorem ~\ref{tomat}, cf.\ \eqref{inflammation}, that
\begin{equation*}
  \int_\U M^\xi \w M^g=
  \int_{\widetilde \U} [\widetilde D_\m] \wedge [\widetilde D]
  \w\frac{1}{1-\varpi^*\pi_F^*\omega_\beta}
  \w\frac{1}{1-\varpi^*\pi_E^*\omega_\alpha}
  \w\frac{1}{1-\varpi^*_\m\pi_{E_\m}^*\omega_\gamma},
\end{equation*}
and, more precisely, for $k<n$,
\begin{equation}\label{mule}
e_k(\M_x)=  \int_\U M^\xi_{n-k} \w M^g_k=
  \int_{\widetilde \U} [\widetilde D_\m] \wedge [\widetilde D]
  \w\frac{1}{1-\varpi^*\pi_F^*\omega_\beta}
  \w\frac{1}{1-\varpi^*\pi_E^*\omega_\alpha}\wedge
  (\varpi^*_\m\pi_{E_\m}^*\omega_\gamma)^{n-k-1};
\end{equation}
here the intersection of divisors should be interpreted in the sense
of Section ~\ref{genz}, so that
\begin{equation*}
  [\widetilde D_\m] \wedge [\widetilde D]
  = [\widetilde D_\m] \wedge
 \1_{\widetilde Y\setminus (p\circ \pi \circ \varpi_\m)^{-1}\{x\}}
  [\widetilde D].
  \end{equation*}
For $k=n$ we have an expression already on $Y$, cf.\
\eqref{stackmoln},
\begin{equation}\label{lyssna}
 e_n(\M_x)=  \int_\U M^\xi_{0} \w M^g_n=
 \int_\U \1_{\{x\}} M^g_n =
 \int_{(p\circ \pi_E)^{-1}\U} \1_{(p\circ \pi_E)^{-1}\{x\}}
[D]\w\frac{1}{1-\pi_F^*\omega_\beta}\w\frac{1}{1-\pi_E^*\omega_\alpha}.
 \end{equation}

As in Section ~\ref{multiplicera}, let $h_E$, $h_F$, and $h_\m$ be the first Segre
classes the tautological line bundles of $\P(E)$, $\P(F)$, and
$\P(E_\m)$,
respectively, and let $\widetilde h_E$, $\widetilde h_F$, and
$\widetilde h_\m$ denote
the pullbacks to $Y$.
Decompose $D$ as in Section
~\ref{segremult},
i.e., $D=D^{X\setminus \{x\}}+D^{\{x\}}$, where
$[D^{X\setminus \{x\}}]=\1_{Y\setminus (p\circ
  \pi_{E})^{-1}\{x\}}  [D]$
and $[D^{\{x\}}]=\1_{(p\circ
  \pi_{E})^{-1}\{x\}}  [D]$, and
decompose $\widetilde D$ analogously.
Then \eqref{lyssna} and \eqref{mule} give the following intersection
formulas, cf.\ Remark \ref{mintgron}.
\begin{prop}\label{kanin}
With the notation above,
\begin{eqnarray*}
e_n(\M_x) &=&
              \int_{(p\circ \pi_E)^{-1}\U}
\sum_{j=0}^{N-1} h_E^j h_F^{N-1-j}  D^{\{x\}},
\\
e_k(\M_x) &=&  
\int_{\widetilde \U} \widetilde h_\m^{n-k-1}\sum_{j=0}^{k+r-2}\widetilde
h_E^j \widetilde h_F^{k+r-2-j} \widetilde D_\m\cdot \widetilde
D^{X\setminus \{x\}}, ~~~~ k<n.
\end{eqnarray*}
\end{prop}
To facilitate the understanding of the second integral, recall that $\dim \widetilde Y=N=n+r-1$.

\begin{ex}\label{taxen}
  Assume that $\M_x$ has support at $\{x\}$.
    Then $D^{\{x\}}=D$ and thus the formula for $e_n(\M_x)$ just
    coincides with the intersection formula \eqref{dimma} for the
    Buchsbaum-Rim multiplicity
    $e(\M_x)$.
\end{ex}

\begin{ex}\label{rodpenna}
Assume that $f$
is a holomorphic section of a Hermitian vector
  bundle $F\to X$ regarded as
    a morphism $g \colon E\to F$, where $E$ is a
    trivial line bundle equipped with the trivial metric.
     We saw in Example \ref{generalisering} that $M^g=M^f$, so we know that $\mult_x M^g_k$ is the
    $k$th Segre number of $\J_x$, where $\J$ is the ideal sheaf
    generated by $f$. In \cite{ASWY} this was proved
    using Vogel cycles.

    Now $\P(E)=X$, $h_E=0$, 
 $\pi_E=\pi \colon Y\to X$ is the blowup of $X$ along $\J$, $h_F=h$ is the
first Segre class of the corresponding tautological line bundle, and
(the top square of) \eqref{stortdiagram} is just
the diagram \eqref{diagram}.
It follows that
Proposition \ref{kanin} just gives back
Gaffney-Gassler's
intersection formulas \eqref{stackmoln} for Segre numbers.
   \end{ex}

In the literature there are various generalizations of the
Buchsbaum-Rim multiplicity. In \cite{CBP07, CBP} Callejas-Bedregal and
Jorge P\'erez introduced 
algebraically defined generalized Buchsbaum-Rim multiplicities
$\tilde e_x(\K_x)$ such that $\tilde e_n(\K_x)$ is the
Buchsbaum-Rim multiplicity of $\M_x:=\Ok_x^s/\K_x$ if $\M_x$ has finite
length, and $\tilde e_x(\K_x)$ is the $k$th Segre number if
$\K_x\subset \Ok_x$ is an ideal.
However there seems to be
no geometric interpretation of these numbers in the
literature. It
would be interesting to compare these generalized Buchsbaum-Rim
multiplicities to the numbers $e_k(\M_x)$. 

\section{The Buchsbaum-Rim multiplicity and Vogel cycles}\label{vogelsektion}

In this section we will construct a generalization of the current
$M^{f,\beta}$ from Section ~\ref{tidigast} and prove a generalization
of \eqref{burra}. In Section ~\ref{advent} we use this to compute $M^g$
in an explicit example.

Assume that we have a morphism $g \colon E \to F$ as in Section
~\ref{sec:MgDef},  and assume
furthermore that $F$ is a trivial bundle, $F=X\times \C^s$.
Let us use the notation from Section
~\ref{sec:MgDef} and let $Z'$ be the zero-set of $G$. Assume that $g$
is nontrivial. Note that then $Z'$ is a proper subvariety
of $\P(E)$,
and moreover that $Z'\subset p^{-1} Z$.
Then $G=(G_1,\ldots, G_s)$ is a section of $p^*F\otimes L_E^{-1}$ that
we can identify with
$L_E^{-1}\oplus \cdots \oplus L_E^{-1}$, and for
$\beta=[\beta_1:\ldots : \beta_s]\in \P^{s-1}$,
$\beta\cdot G:= \beta_1 G_1+\cdots +\beta_s G_s$
is a (homogeneous in $\beta$) section of $L_E^{-1}$. Let $a$ be a nonvanishing
holomorphic section of $L_E$. Then $\beta \cdot Ga$ is a (homogeneous
in $\beta$)
function defining the divisor defined by $\beta\cdot G$, and thus
$dd^c \log |\beta \cdot Ga|^2=[\beta \cdot G]$, where $|\cdot |$ now
denotes absolute value; note that this is
independent of the homogeneous coordinates $\beta$ on $\P^{s-1}$.
For $\beta=(\beta_1,\ldots, \beta_N) \in (\P^{s-1})^N$
and $k=1,\dots,N$,
let
\begin{equation*}
  \Mboll^{G,\beta}_k = \1_{Z'} [\beta_k \cdot G] \wedge \dots
  \wedge [\beta_1 \cdot G],
   \end{equation*}
 cf.\ Section ~\ref{tidigast},
and let
$\Mboll^{G,\beta} = \sum_{k=1}^N \Mboll^{G,\beta}_k$. Moreover, in analogy
with \eqref{Mgdef} define
\begin{equation}\label{himmelsbla}
    M^{g,\beta} = p_*\left(\Mboll^{G,\beta}\w\frac{1}{1-\omega_\alpha} \right),
  \end{equation}
and for $k=0, \ldots, n$, let $M^{g,\beta}_k$ denote the component of bidegree $(k,k)$.

If $g\equiv 0$, then $G\equiv 0$, and thus we can define
$\Mboll^{G,\beta}=\Mboll^{G,\beta}_0=\1_{\P(E)}$. Defining $M^{g,\beta}$
through \eqref{himmelsbla}, we get
\begin{equation}\label{lukta}
  M^{g,\beta}=p_*\left ( \frac{1}{1-\omega_\alpha}\right )=M^g,
\end{equation}
see Example
\ref{triviala}.

  \begin{ex}
Assume that $E$ is a trivial line bundle, cf.\ Example
~\ref{generalisering}. Then $g$ can be identified with an $s$-tuple
$f$ of holomorphic functions, and one can check that
\begin{equation*}
  M^{g,\beta}_k =  \Mboll^{G,\beta}_k = M^{f,\beta}_k=
  \1_Z  [\beta_k \cdot f] \wedge \dots
   \wedge [\beta_1 \cdot f],
 \end{equation*}
 where $Z$ is the zero-set of $f$.
\end{ex}

First, we have the following generalization of \eqref{hurra}.
\begin{prop} \label{maska}
Assume that $g \colon E \to F=X\times \C^s$ is pointwise injective
outside $\{x\}$.
Then, for a generic $\beta \in
(\P^{s-1})^N$,
\begin{equation}\label{linjer}
  M^{g, \beta} = M^g.
\end{equation}
\end{prop}

  The proof of Proposition \ref{maska}
    is inspired by the proof of
  Theorem~6.1 in \cite{ASWY}.

\begin{proof}
We will show that for a generic $\beta\in (\P^{s-1})^N$ and for
$k=1, \dots,N$
\begin{equation} \label{Mbollddc}
 \Mboll_k^{G, \beta} -  \Mboll_k^G
   = dd^c w_k',
\end{equation}
where $w_k'$ is a $(k-1,k-1)$-current that is the push-forward of an
intrinsic current on $Z'$.
It follows that
\begin{equation*}
  M^{g,\beta}-M^g=
  M^{g,\beta}_n - M^g_n = dd^c w_n,
\end{equation*}
where $w_n$ is a $(n-1,n-1)$-current that is the push-forward of an
intrinsic current on $Z=\{x\}$. Now $w_n$ vanishes for degree reasons
and thus \eqref{linjer} holds.

Let $\J$ be the ideal sheaf generated by $G$ and let $\pi \colon Y\to \P(E)$ be a
principalization of $\J$ with exceptional
divisor $D$.
Then $\pi^*G=G^0G'$, where $G^0$ is a section of the line bundle $\L \to
Y$ that defines $D$ and $G'$ is a nonvanishing section of
$\pi^*(p^*F\otimes L_E^{-1})\otimes \L^{-1}$.
Let $a$ be a (local) nonvanishing section of $L_E$.
Then
\begin{equation}\label{orange}
  dd^c\log|\pi^*(Ga)|^2=[D]+\omega_G,
\end{equation}
where $\omega_G$ is as in Remark \ref{borste}.
Hence
\begin{equation}\label{apelsin}
  \Mboll^G = \pi_*\left ([D]\w \frac{1}{1-\omega_G}\right ).
  \end{equation}

Since $Z=\{x\}$, $Z'$ is compact, and so is $|D|$. In particular, $D$ consists of
finitely many irreducible components. It follows that for a generic
$\beta=(\beta_1,\ldots,
\beta_N) \in (\P^{s-1})^N$, $D$ and
the divisors of the sections $\beta_k \cdot G',\dots,\beta_1 \cdot G'$
of $\pi^*L_E^{-1}\otimes \L^{-1}$
intersect properly for each $k$. Fix such a $\beta$.
In view of \eqref{brajobbat} and \eqref{pushpull}, for $k=1,\ldots, N$,
\begin{equation}\label{begynnelsen}
  [\beta_k\cdot G]\wedge \cdots \wedge [\beta_1\cdot G]
  =
 \pi_* \left ([D]\w [\beta_{k-1}\cdot G']\wedge \cdots \wedge [\beta_1\cdot G']
  +
   [\beta_k\cdot G']\wedge \cdots \wedge [\beta_1\cdot G']\right ),
\end{equation}
  and thus
  \begin{equation}\label{ordentligt}
  \Mboll^{G,\beta}_k= \pi_* \left ([D]\w [\beta_{k-1}\cdot G']\wedge \cdots
  \wedge [\beta_1\cdot G'] \right ),
\end{equation}
cf.\ \eqref{ordobild}.
  In particular $\Mboll^{G, \beta}_1=\pi_*[D] =  \Mboll^G_1$,
  cf.\ \eqref{apelsin},
  so that \eqref{Mbollddc} is trivially satisfied for $k=1$.

For a local nonvanishing section $a$ of $L_E$ and $\ell=1,\ldots, N$, let
$$
\theta_\ell 
=\pi^*\log (|\beta_\ell\cdot Ga|^2/|Ga|^2).
$$
As in Section \ref{mf}, the metrics on $E$ and $F$ induce metrics on $\pi^*L_E^{-1}\otimes \L^{-1}$ and
$\pi^*(p^*F \otimes L_E^{-1})\otimes \L^{-1}$ such that
\[
  \theta_\ell = \log (|\beta_\ell\cdot G'|^2/|G'|^2).
  \]
In particular, $\theta_\ell$ is a globally defined function on $Y$.
Next, for $k=2,\ldots, N$, let
\[
w_k'' := \sum_{\ell=1}^{k-1} \theta_\ell [D]  \wedge [\beta_{\ell-1}\cdot G'] \wedge \cdots \wedge [\beta_1 \cdot G']\wedge \omega_G^{k-1-\ell};
\]
for $\ell=1$ the summand is interpreted as
$\theta_1[D]\w\omega_G^{k-2}$.
Then $w_k''$ is
a $(k-1,k-1)$-current on $Y$ that is the pushforward of an intrinsic current
on $D$, cf.\ \eqref{middag}.
By the Poincar\'e-Lelong formula, cf.\ \eqref{gamas1},
\begin{equation*}
  dd^c \log |\pi^*(\beta_\ell \cdot Ga)|^2=
  dd^c \log |G^0\beta_\ell \cdot G'\pi^*a|^2
  =[D]+[\beta_\ell\cdot G'],
\end{equation*}
and hence
\[
  dd^c w_k''=
  [D]\w [\beta_{k-1}\cdot G']\wedge \cdots \wedge [\beta_1\cdot G']-
  [D]\w \omega_G^{k-1},
\]
cf.\ \eqref{orange}.
In view of \eqref{apelsin} and \eqref{ordentligt}, thus
$w_k':=\pi_*w_k''$ satisfies \eqref{Mbollddc} and is of the desired
form, since $D$ is mapped onto $Z'$.
\end{proof}

For general $g$ we have the following generalization of
\eqref{burra}.
\begin{prop}\label{slaska}
  For a generic $\beta\in (\P^{s-1})^N$ and $k=0,\ldots, n$,
  \begin{equation}\label{braska}
\mult_x M^{g,\beta}_k = \mult_x M^{g}_k.
    \end{equation}
  \end{prop}

  \begin{proof}
Note that if $g\equiv 0$, then the proposition is trivially true
in view of \eqref{lukta}. We may therefore assume that $g$ is
nontrivial.

    Let us first consider the case $k=n$. Then it suffices to prove
    that, for a generic $\beta\in (\P^{s-1}) ^N$ and $\ell=1,\ldots, N$,
    \begin{equation}\label{tennisberg}
 \1_{p^{-1}\{x\}} \Mboll_\ell^{G, \beta} -  \1_{p^{-1}\{x\}} \Mboll_\ell^G
   = dd^c v_\ell',
\end{equation}
where $v_\ell'$ is a $(\ell-1,\ell-1)$-current that is the push-forward of an
intrinsic current on $p^{-1}\{x\}$.
Then it follows as in the proof of Proposition \ref{maska} that
 $\1_{\{x\}}M^{g,\beta}_n - \1_{\{x\}}M^g_n = 0$,
 and consequently $\mult_x M^{g,\beta}_n = \mult_x M^{g}_n$, cf.\
 \eqref{lyssna}.
To prove \eqref{tennisberg} consider the setting
 in the proof of Proposition ~\ref{maska}. Since
    $(p\circ \pi) ^{-1}\{x\}$ is compact, if we replace $X$ by a small
    enough
    neighborhood of $x$, then
    for a generic $\beta\in (\P^{s-1})^N$, $\beta_{\ell}\cdot G',\ldots,
    \beta_1\cdot G'$ and $D$ intersect properly 
    for each
    $\ell$. For such a $\beta$ let $v_\ell''=\1_{(p\circ \pi)^{-1}\{x\}} w_\ell''$, where $w_\ell''$ is as
    in the proof of  Proposition \ref{maska}. Then by arguments as in
    that proof $v_\ell':=\pi_* v_\ell''$ is of the desired form.

    For $k<n$, let us consider the setting of Section ~\ref{hundra}, and let us
    use the shorthand notation $\widetilde \omega_\alpha, \widetilde
    \omega_\beta, \widetilde \omega_\gamma$ for the pullbacks to
    $\widetilde Y$ of the
    forms  $\omega_\alpha, \omega_\beta, \omega_\gamma$.
    We can choose the principalization of $G$ in the proof of
    Proposition ~\ref{maska} to be $\pi:=\pi_E\circ \varpi: \widetilde
    Y\to \P(E)$, i.e, $Y$ in that proof translates to $\widetilde
    Y$; similarly $D$ translates to $\widetilde D$ and $\omega_G$ to
    $\widetilde \omega_\beta$.

    First, assume that $r\ge 2$. Let $\theta_\ell$ be as in the proof
    of Proposition ~\ref{maska} and,
    for $k=0,\ldots, n-1$,
    let
\begin{equation}\label{sommardag}
  w_k'' = \sum_{j=1}^{k+r-2}\sum_{\ell=0}^{j-1}
  \theta_{\ell+1} [\widetilde D_\m]\w[\widetilde D]
  \wedge [\beta_{\ell}\cdot G'] \wedge \cdots \wedge [\beta_1
  \cdot G']
  \wedge
  \widetilde\omega_\beta^{j-1-\ell}\w
  \widetilde\omega_\alpha^{k+r-2-j}\w
  \widetilde\omega_\gamma^{n-k-1};
\end{equation}
for $\ell=0$ the summands are interpreted as
$\theta_1 [\widetilde D_\m]\w[\widetilde D]
  \wedge
  \widetilde\omega_\beta^{j-1}\w
  \widetilde\omega_\alpha^{k+r-2-j}\w
  \widetilde\omega_\gamma^{n-k-1}.
  $
Now, as above, for a generic $\beta\in (\P^{s-1})^N$,
$\beta_k\cdot G',\ldots, \beta_1\cdot G'$ and $[\widetilde
D_\m]\w[\widetilde D]$ intersect properly for each $k$,
and thus
\begin{multline*}
  dd^c w_k''
  =
  \sum_{j=1}^{k+r-2}
   [\widetilde D_\m]\w[\widetilde D]
  \wedge [\beta_{j}\cdot G'] \wedge \cdots \wedge [\beta_1
  \cdot G']
  \w  \widetilde\omega_\alpha^{k+r-2-j}\w
  \widetilde\omega_\gamma^{n-k-1}
  \\
-  \sum_{j=1}^{k+r-2}
  [\widetilde D_\m]\w[\widetilde D]
  \wedge
  \widetilde\omega_\beta^{j}\w
  \widetilde\omega_\alpha^{k+r-2-j}\w
  \widetilde\omega_\gamma^{n-k-1}.
\end{multline*}
 Using \eqref{lurig}, \eqref{inflammation}, Theorem \ref{tomat},
 \eqref{himmelsbla}, and \eqref{begynnelsen} one can check that
\[
  dd^c p_*\pi_* w_k'' = M^\xi_{n-k}\w M^{g,\beta}_k -
  M^\xi_{n-k}\w M^{g,}_k,
\]
cf.\ \eqref{mule}.
Now for each $k$, $w_k''$ is an $(N-1, N-1)$-current on $\widetilde Y$
that is the pushforward of an intrinsic current on
$\widetilde D$, cf.\ \eqref{middag}, and thus $p_*\pi_* w_k''$ is the pushforward on an
intrinsic $(n-1, n-1)$ current on $\{x\}$ and therefore it vanishes
for degree reasons.
Thus $M^\xi_{n-k}\w M^{g,\beta}_k = M^\xi_{n-k}\w M^{g,}_k$, which
by \eqref{ballong} implies \eqref{braska}.

Next, assume that $r=1$. Then the assumption $g\not\equiv
      0$ implies that $g$ is generically injective, i.e., $Z$ is a
      proper subvariety of $X$ and thus $M^g_0=M^{g,\beta}_0=\1_Z=0$,
      so that \eqref{braska} is trivially satisfied. For $k=1,\ldots,
      n-1$ the proof above works verbatim. In fact, if
      $r+k=2$, the left sum in \eqref{sommardag} is empty, and thus
      $w_k''=0$. In this case $M_k^g$ coincides with
      $M^{g,\beta}_k$ since they are both the pushforward of
      $[\widetilde D]$.

    \end{proof}

\section{An explicit example}\label{advent}

Let $\Ok_0$ denote the ring of germs of holomorphic functions on $\C^2_{(x_1,x_2)}$ at $0$,
and let $\M_0=\Ok_0^2/\Im g^*$ where $g^*\colon \Ok_0^3\to \Ok_0^2$ is given
by
\begin{equation}\label{matrix}
g^\ast = \begin{bmatrix}
	x^a_1 & 0 & x_2^c \\
	0 & x_2^b  &   x_1^d\\
\end{bmatrix},
\end{equation}
where $a,b,c,d\ge 1$.
Now $g^*$ is a so-called parameter matrix, i.e.
an $r \times (r+n-1)$-matrix, that is pointwise surjective
outside of $0$, and therefore
$e(\M_0) = \dim_\C \M_0$, see e.g.,
\cite[Corollary~4.5]{BR2}.
By calculating $\dim_\C \M_0$ one can verify that
\begin{equation}\label{snut}
e(\M_0) = ab+ac +bd.
\end{equation}

\smallskip

We will compute $e(\M_0)$ using Theorem \ref{thm1} and Proposition ~\ref{maska}.
Let $X=\C^2_{(x_1,x_2)}$ and let
$E=X\times \C^2_{(\alpha_1,\alpha_2)}$ and $F=X\times \C^3$, equipped with the trivial metrics.
Then \eqref{matrix} gives a holomorphic morphism $g \colon E\to F$
such that
$\M_0= \Ok(E^*)_0/\Im g^*_0$.

By Theorem ~\ref{thm1} and Proposition ~\ref{maska}
it is enough to compute
$M^{g,\beta}$ for a generic choice of $\beta\in (\P^2)^3$. We will start by computing
$\Mboll^{G,\beta}$.
Now
\begin{equation}
  G(x,[\alpha])=g\alpha =
  (\alpha_1x_1^a, \alpha_2x_2^b, \alpha_1x_2^c+\alpha_2x_1^d)=:
(G_1, G_2, G_3),
\end{equation}
where  each $G_j$ is a function that is $1$-homogeneous
  in $\alpha$, i.e., a section of $L_E^{-1}$.
  Note that
  \begin{equation}\label{davidsson}
    Z'=\{x=0\}=\{0\}\times \P^1\subset X\times
    \P^1=\P(E)
  \end{equation}
  has codimension $2$.
Thus
$\Mboll^{G,\beta}= \Mboll^{G,\beta}_2+\Mboll^{G,\beta}_3$.

Note that from a generic $\beta'=(\beta'_1,\beta'_2,\beta'_3) \in  ( \mathbb{P} ^{2})^3$, by
using row operations
\begin{equation*}
  \beta'_\ell\cdot G+c\beta'_k\cdot G, ~~~~~ k<\ell, ~~~~~ c\in \C
  \end{equation*}
we can obtain $\beta
\in  ( \mathbb{P} ^{2})^3$ of the form
\begin{equation}\label{balkong}
  \beta =(\beta_1,\beta_2,\beta_3)= \big ([\beta_{11}:\beta_{12}:1], [\beta_{21}:1:0],
  [1:0:0]\big ).
  \end{equation}
In view of \eqref{lussekatt},
$\Mboll^{G,\beta}_k=\Mboll^{G,\beta'}_k$. We may therefore assume
that $\beta$ is of the form \eqref{balkong} and moreover, since
$\beta'$ is generic, that
$\beta_{21}\neq 0$.
Let $h= (h_1,h_2,h_3)=(\beta_1 \cdot G, \beta_2\cdot G, \beta_3\cdot G)$, i.e.
\begin{align*}
    h_1 &= \beta_{11} \alpha_1x_1^a + \beta_{12} \alpha_2 x_2^b + \alpha_1x_2^c + \alpha_2x_1^d \\
h_2 &= \beta_{21} \alpha_1x_1^a + \alpha_2 x_2^b \\
h_3 &= \alpha_1x_1^a.
\end{align*}
Since $Z'$ is compact, for a generic $\beta\in (\P^2)^3$, $h$
is a Vogel sequence
for the ideal sheaf $\J$ generated by $G$ at each point $(0,[\alpha])\in
Z'$, cf.\ Section \ref{tidigast}.
In particular,
the intersection $[h_2]\wedge [h_1]$ is proper and thus of the form
\begin{equation} \label{propra}
    [h_2]\wedge [h_1] = m [Z'] + [\mathcal{C}],
  \end{equation}
  where $\mathcal C$ is a cycle of codimension $2$ not contained in
  $Z'$.
  It follows that
  \begin{equation}\label{lilla}
    \Mboll^{G,\beta}_2=\1_{Z'} [h_2]\wedge [h_1] = m [Z'].
  \end{equation}

  Next, we consider
  \begin{equation}\label{stora}
  \Mboll^{G,\beta}_3= \1_{Z'} [h_3] \wedge [h_2]\wedge [h_1]  =
  \1_{Z'}[h_3]\wedge [\mathcal{C}]  = \1_{Z'} [\alpha_1]\wedge
  [\mathcal{C}] + a\1_{Z'} [x_1] \wedge [\mathcal{C}].
    \end{equation}
First, note that $[\alpha_1]$ intersects $[h_2]\wedge [h_1]$ at the point
$\mathbf 0:= (0,0,[0:1])$ (independently of the $\beta_{ij}$). In particular, the
intersection is proper.
By \eqref{eq:MArestriction} we get
\begin{multline*}
  [\alpha_1]\wedge[h_2]\wedge [h_1]
    =
 [\alpha_1]\w [h_2(x_1, x_2, [0:1] )] \wedge [h_1(x_1, x_2, [0:1])]
  =\\
 [\alpha_1]\w[x_2^b]\wedge [\beta_{12} x_2^b + x_1^d]
 =
[\alpha_1]\w[x_2^b]\wedge [x_1^d]
  = bd [\mathbf 0].
\end{multline*}
Hence, by \eqref{propra} and \eqref{davidsson},
\begin{equation} \label{hosianna}
 \1_{Z'}[\alpha_1]\wedge [\mathcal{C}] =\1_{Z'} [\alpha_1]\wedge[h_2]\wedge [h_1] - m \1_{Z'}[\alpha_1]\wedge[Z']  = (bd - m)  [\mathbf 0].
\end{equation}
Next, since $h$
is a Vogel sequence for $\J$ at each point in $Z'$ and $Z'$ is compact,
 $[x_1]\wedge [\mathcal C]$ has support at a finite number of points in $Z'$.
It follows from, e.g., \cite[p.\ 137]{Ch}, that for a generic small $
\epsilon >0$ the number of points counted with multiplicities equals the
number of intersection points of (the proper intersection) of
$[x_1-\epsilon]$ and $ [\mathcal{C}]$.
Since $[x_1-\epsilon]$ does not intersect $[Z']$,
  $[x_1-\epsilon]\wedge  [\mathcal{C}]$
  equals the proper intersection
  $[x_1-\epsilon]\wedge  [h_2]\wedge [h_1].$
Since, by assumption, $\beta_{21}\neq 0$,
\[
\{h_2(\epsilon, x_2, [\alpha])=0\} = \{[\alpha]=
[-x_2^b/(\beta_{21}\epsilon^a): 1]\}.
\]
Thus from
\eqref{eq:MArestriction} we get
\begin{multline*}
  [x_1-\epsilon] \wedge [h_2]\wedge [h_1]=
  [x_1-\epsilon] \wedge [h_2(\epsilon, x_2, [\alpha])]\wedge
  [h_1(\epsilon, x_2, [\alpha])]
  =\\
  [x_1-\epsilon] \wedge [h_2(\epsilon, x_2, [\alpha])]\wedge
  [h_1(\epsilon, x_2, [-x_2^b/(\beta_{21}\epsilon^a):1]) ].
  \end{multline*}
Since
\begin{equation*}
 0=   h_1(\epsilon, x_2,  [-x_2^b/(\beta_{21}\epsilon^a):1]) = \left(\beta_{12} - \frac{\beta_{11}}{\beta_{21}} \right)x_2^b - \frac{1}{\beta_{21} \e^a} x_2^{b+c} + \e^d
\end{equation*}
has $b+c$ solutions, we conclude that $[x_1-\e]$ and
$[\mathcal C]$ intersect at $b+c$ points.
Hence
\begin{equation} \label{hatet}
\int_{\P(E)}\1_{Z'} [x_1] \wedge [\mathcal{C}] = \int_{\P(E)} [x_1] \wedge [\mathcal{C}] = \int_{\P(E)} [x_1-\epsilon] \wedge [\mathcal{C}] = b+c
\end{equation}
and by combining \eqref{stora}, \eqref{hosianna}, and
\eqref{hatet} we get
\begin{equation}\label{hemmabio}
  \int_{\P(E)} \Mboll^{G,\beta}_3=
  \int_{\P(E)}
  \1_{Z'} [\alpha_1]\wedge [\mathcal{C}]
  +
  a \int_{\P(E)} \1_{Z'} [x_1] \wedge [\mathcal{C}]
  =
  bd-m+ a(b+c).
    \end{equation}

We conclude
 \begin{multline*}
   e(\M_0)=
      \int_X M^{g} =\int_X M^{g,\beta}=
      \int_{\P(E)}
   \Mboll^{G,\beta} \w \frac{1}{1-\omega_\alpha} =\\
   \int_{\P (E)} \Mboll^{G,\beta}_2 \w \omega_\alpha+
   \int_{\P (E)} \Mboll^{G,\beta}_3
   =  m + bd-m +  a(b+c)= ab + ac +bd.
 \end{multline*}
Here we have used Theorem ~\ref{thm1} for the first equality,
Proposition ~\ref{maska}
for the second equality,
\eqref{himmelsbla} for the third equality, and \eqref{lilla} and
\eqref{hemmabio} for the last equality.
Hence we get \eqref{snut} and
$M^g = (ab + ac +bd)[0]$.

\def\listing#1#2#3{{\sc #1}:\ {\it #2},\ #3.}

\end{document}